\documentclass[10pt,oneside,leqno]{amsart}
\usepackage{amsxtra}
\usepackage{amsopn}
\usepackage{amsmath,amsthm,amssymb}
\usepackage{amscd}
\usepackage{amsfonts}
\usepackage{latexsym}
\usepackage{verbatim}

\theoremstyle{plain}
\newtheorem{theorem}{Theorem}[section]

\newtheorem{prop}[theorem]{Proposition}

\newtheorem{rem}[theorem]{Remark}
\newtheorem{ex}[theorem]{Example}

\newcommand\R{{\mathbb R}}
\newcommand\Z{{\mathbb Z}}
\renewcommand\H{{\mathbb H}}
\newcommand\T{{\mathbb T}}
\begin{document}
\title{Generalized $G_2$-manifolds  and $SU(3)$-structures}
\author{Anna Fino and Adriano Tomassini }
\date{\today}
\address{Dipartimento di Matematica \\ Universit\`a di Torino\\
Via Carlo Alberto 10\\
10123 Torino\\ Italy} \email{annamaria.fino@unito.it}
\address{Dipartimento di Matematica\\ Universit\`a di Parma\\ Viale G.P. Usberti 53/A\\
43100 Parma\\ Italy} \email{adriano.tomassini@unipr.it}
\thanks{This work was supported by the Projects MIUR ``Geometric Properties of Real and Complex Manifolds'',  "Riemannian Metrics  and Differentiable
 Manifolds" and by GNSAGA
of INdAM}
\begin{abstract}
We construct a family of compact   $7$-dimensional manifolds
endowed with a weakly integrable generalized $G_2$-structure with
respect to a closed and non-zero $3$-form. We relate the previous
structures with $SU(3)$-structures in dimension 7. Moreover, we
investigate which types of $SU(3)$-structures on a $6$-dimensional
manifold $N$ give rise to a strongly integrable generalized
$G_2$-structure with respect to a non-zero $3$-form on the product
 $N \times S^1$.

\end{abstract}
\maketitle
\section{Introduction}

The notion of generalized geometry goes back to the work of Hitchin \cite{Hi1} (see also \cite{Hi3}).  In this context,
 Witt  \cite{W} introduced a  new type of structures on a $7$-dimensional
manifold  $M$ in terms of a differential form of mixed degree,
thus generalizing the classical notion of $G_2$-structure
determined by a stable and positive $3$-form. Instead of studying
geometry on the tangent bundle $TM$ of the manifold,  one
considers the   bundle   $TM \oplus T^*M$ endowed with a natural
orientation and an  inner product  of signature $(7, 7)$,  where
$T^*M$ denotes  the cotangent bundle of $M$. In this way,  if $M$
is spin,  then the differential form of mixed type can be viewed
as a $G_2 \times G_2$-invariant spinor  $\rho$ for the  bundle and
it is  called  the structure form.

These  structures are called {\it generalized $G_2$-structures}
and they induce a Riemannian metric, a $2$-form $b$ (the
$B$-field), two unit spinors $\Psi_{\pm}$   and a function $\phi$
(the dilaton).   By \cite{W},  any  $G_2 \times G_2$-invariant
spinor   $\rho$ is stable and has a canonical  expression by $\rho
= e^{-\phi} e^{\frac{b}{2}} \wedge (\Psi_+ \otimes
\Psi_-)^{ev,od}$ in terms of the two spinors, the $B$-field and
the dilaton function. In the paper we will  restrict to  the case
of constant dilaton, i.e. $\phi = {\mbox const}$, and trivial
$B$-field.

Up to a $B$-field transformation, a generalized $G_2$-structure is
essentially a pair of $G_2$-structures.  If the two spinors
$\Psi_+$ and $\Psi_-$ are linearly independent, then the
intersection of  the two isotropy groups, both  isomorphic to
$G_2$, determined  by the two spinors coincides with $SU(3)$.
Therefore,  one can express the structure form in terms of  the
form $\alpha$ dual to the unit vector stabilized by $SU(3)$ and of
the forms $(\omega, \psi = \psi_+ + i \psi_-)$, associated  with
$SU(3)$, where $\omega$ is the fundamental form and  $\psi$  is
the complex volume form. Assuming that the angle between $\Psi_+$
and $\Psi_-$ is $\frac {\pi} {2}$, then it turns out that
\begin{equation} \label{expression}
\begin{array}{l}
\rho = (\Psi_+ \otimes \Psi_-)^{ev} = \omega +  \psi_+ \wedge
\alpha -  \frac{1}{6}  \omega^3 \wedge \alpha,\\ [5pt] \hat \rho =
(\Psi_+ \otimes \Psi_-)^{od} = \alpha -  \psi_-  -  \frac{1}{2}
\omega^2 \wedge \alpha,
\end{array}
\end{equation}
where $\hat \rho$ is the companion of $\rho$ and $\omega^k$
denotes the $k$-power wedge of $\omega$. In this paper we will
consider generalized $G_2$-structures defined by the previous
structure forms.  In this case, the two associated
$G_2$-structures do not coincide.

If $H$ is a $3$-form  (not necessarily closed) on $M$, then one can
 consider two
types of  generalized $G_2$-structures with respect to the
$3$-form $H$: \newline the {\em strongly integrable} ones, i.e.
those associated to  a structure form $\rho$  which satisfies
\begin{equation} \label{strongnotclosed}
d  \rho  + H \wedge \rho =  d \hat \rho + H \wedge \hat \rho  =0,
\end{equation}
 and the {\em
weakly integrable} ones, i.e. those defined by the condition $$ d
\rho  + H \wedge \rho = \lambda \hat \rho, $$ where $\lambda$ is a
non-zero constant. The  previous structures are said  of even or
odd type according to the parity of  $\rho$.

Note that these definitions of integrability are slightly
different from the ones given  in \cite{W}, where  the closure of
the $3$-form $H$ is assumed.

If $H$ is closed, then   the twisted operator $d_H   \cdot =  d
\cdot + H \wedge \cdot$ defines a differential complex and if, in
addition, $M$ is compact, then  the strongly integrable
generalized $G_2$-structures can be interpreted  as critical
points of a  a certain  functional  \cite[Theorem 4.1]{W}. In this
case  the underlying spinors $\Psi_{\pm}$ are parallel with
respect to the Levi-Civita connection and therefore there exist no
non-trivial compact examples  with such structures, i.e. there are
only the classical examples of manifolds with   holonomy contained
in $G_2$. If $H$ is not closed, then we will show that compact
examples can be constructed starting from a $6$-dimensional
manifold endowed with an $SU(3)$-structure.

If $H$ is closed, then the weakly integrable generalized
$G_2$-structures can be also viewed as critical points of a
functional under a constraint, but they have no classical
counterpart.  The existence of weakly integrable generalized
$G_2$-structures  with respect to a closed $3$-form $H$ on a
compact manifold  was  posed as    an open problem in \cite{W}. We
 construct such  structures on a family of compact manifolds
and we  relate them with $SU(3)$-structures in dimension $7$,
where $SU(3)$ is identified with the subgroup $SU(3) \times \{ 1
\}$ of $SO(7)$.

After reviewing the general theory of generalized $G_2$-structures, in
 section 3 we  construct  a family of  compact  $7$-dimensional manifolds endowed
 with a weakly integrable generalized $G_2$-structure with respect to
 a closed and non-zero $3$-form $H$ (Theorem \ref{example}).
 The corresponding structure form is the odd type form  $\hat \rho$  given by  \eqref{expression}. These manifolds are
 obtained as a compact quotients $M_{\beta}$ by  uniform discrete subgroups (parametrized by  the  p-th roots of  unity $e^{i \beta}$) of
 a semi-direct product $SU(2) \ltimes
 \H$, where  $\H$
 denotes the quaternions. It turns out that these manifolds have an $SU(3)$-structure $(\omega, \eta, \psi)$ such that
\begin{equation}\label{hyposystem}
d \eta = \lambda \omega, \quad d (\eta \wedge \psi_{\pm}) =0.
\end{equation}
In particular they are contact metric. The structures satisfying
the condition \eqref{hyposystem} can arise on hypersurfaces of
$8$-dimensional manifolds with an integrable $SU(4)$-structure and
they are the analogous of the \lq \lq hypo\rq \rq
$SU(2)$-structures in dimension $5$ (see \cite{CS2}). In the same
vein of \cite{Hi1}, we consider a family $(\omega(t), \eta(t),
\psi(t))$ of $SU(3)$-structures containing the $SU(3)$-structure
$(\omega, \eta, \psi)$ and the corresponding evolution equations.
In this way  in section 4 we
 show that on the product of $M_{\beta}$ with an open interval
there exists a Riemannian metric with  discrete holonomy
contained  in $SU(4)$ (Theorem 4.1).

Starting from a $6$-dimensional manifold $N$ endowed with an
$SU(3)$-structure $(\omega, g, \psi)$, it is possible to define in
a natural way a generalized $G_2$-structure with  the structure
form $\rho$ of even type given by \eqref{expression} on the
Riemannian product $(M = N \times S^1, h)$, with
$$
h = g + dt \otimes dt
$$
and $\alpha = dt$.
 In \cite{W} an example of this type  with
a $6$-dimensional nilmanifold $N$ was considered in order to
construct a compact manifold endowed  with a  strongly integrable
generalized $G_2$-structure with respect to a non-closed $3$-form
$H$.

We will prove in general  that if $N$  is a $6$-dimensional
manifold endowed with an $SU(3)$-structure $(\omega, g, \psi)$,
then the generalized $G_2$-structure defined by $\rho$ on $N
\times S^1$ satisfies the conditions \eqref{strongnotclosed}, for
a non-zero $3$-form $H$, if and only if
\begin{equation}\label{symplectichalfflat}
d \omega =0,  \quad  d \psi_+  = - \pi_2 \wedge \omega, \quad d \psi_- =0,
\end{equation}
where the $2$-form $\pi_2$ is the unique  non zero  component of
the intrinsic torsion (see Theorem \ref{stronglyint}). We will
call $SU(3)$-structures which satisfy the previous conditions
belonging to the class ${\mathcal W}_2^+$.  The $3$-form $H$ is
related to
 the component $\pi_2$ of the intrinsic torsion by $ H = \pi_2
 \wedge \alpha$ and we
will show that $H$ will never be closed unless $\pi_2 = 0$.

It has to be noted that, if $(\omega, g,  \psi)$ is in the class
${\mathcal W}_2^+$, then  the $SU(3)$-structure given by $(\omega,
g,  i \psi)$ is symplectic half-flat (see \cite{CS}), i.e. the
fundamental form $\omega$ and the real part of the  complex volume
form are both closed.  The half-flat structures turn out to be
useful in the construction of metrics with holonomy group
contained in $G_2$ (see e.g. \cite{Hi1,CS,CF}). Indeed, starting
with a half-flat structure on $N$, if certain evolution equations
are satisfied, then there exists a Riemannian metric with holonomy
contained in $G_2$ on the product of the manifold $N$ with some
open interval. Examples of compact manifolds with symplectic
half-flat structures have been given in \cite{CT}, where invariant
symplectic half-flat structures on nilmanifolds are classified.
Other examples are considered in \cite{dBT}  where Lagrangian
submanifolds are studied instead.

\medskip

{\em{Acknowledgements}}
The authors thank Simon Salamon, Frederik Witt   for useful comments and suggestions. They also thank the  Centro di Ricerca Matematica \lq \lq Ennio De Giorgi\rq \rq, Pisa,  for the warm hospitality.

\section{Generalized $G_2$ structures and spinors}

In this section we are going to recall some facts on generalized
$G_2$-structures which have been studied by Jeschek  and Witt in
\cite{W,W2,JW} in the general case of $\phi$ non-constant and
non-trivial $B$-field. In the next sections we will deal with the
case $\phi = const$ and trivial $B$-field.

 Let $V$ be a
$7$-dimensional real vector space and denote by $V^*$  the dual
space of $V$. Then $V \oplus V^*$ has a natural  orientation and a
inner product of signature $(7,7)$ defined by
$$
(v + \xi, v + \xi) = -  \frac 12 \xi(v), \quad \forall v \in V, \,
\xi \in V^*.
$$
The inner product determines a group coniugate to $SO(7,7)$ inside the linear group
$GL(14)$. Since as $GL(7)$-space ${\mathfrak {so}} (7,7) = End(V)
\oplus \Lambda^2 V^* \oplus \Lambda^2 V$, any $b \in \Lambda^2
V^*$ defines an element (called {\em B-field}) in ${\mathfrak
{so}} (7,7)$. By exponentiating to $SO(7,7)$ the action of
$\Lambda^2 V^* \subset {\mathfrak {so}} (7,7)$
$$
v  \to    v \lrcorner  b,
$$
 one gets an action on $V \oplus V^*$, given by $ \exp (b) (v \oplus
 \xi) = v \oplus
 ( v  \lrcorner  b+ \xi)$.
 Then  $V \oplus V^*$ acts  on $\Lambda^* V^*$ by
$$
(v + \xi) \eta = \iota (v) \eta + \xi \wedge \eta,
$$
and we have
$$
(v + \xi)^2 \eta = - (v + \xi, v + \xi) \eta.
$$
 Therefore $\Lambda^* V^*$ can be viewed  as  a module over the  Clifford algebra of $V \oplus V^*$. The space
  $\Lambda^*  V^*$,   as the spin
 representation of  $Spin(7, 7)$,    determines the splitting of  $\Lambda^* V^* \otimes (\Lambda^7 V)^{\frac 12}$
 $$
 \begin{array} {l}
 S^+ = \Lambda^{ev} V^* \otimes (\Lambda^7 V)^{\frac 12}\\
 S^- = \Lambda^{od} V^* \otimes (\Lambda^7 V)^{\frac 12}
 \end{array}
 $$
 into the sum of the two irreducible spin representations.
By considering $b \in \Lambda^2 V^*$, then one has the following
induced action on
 spinors given by
 $$
 \exp (b) \eta = (1 + b + \frac 12 b \wedge b + \cdots) \wedge \eta = e^b \wedge \eta.
 $$
   If $\sigma$ is the Clifford algebra anti-automorphism defined by $\sigma (\gamma^p) = \epsilon(p) \gamma^p$, on any element of degree $p$, with
 $$
  \epsilon  (p) = \left\{  \begin{array} {cll}
 1 \quad &{\mbox{for}} &\quad p \equiv 0, 3 \quad {\mbox{mod}} \, 4,\\
 - 1 \quad &{\mbox{for}}& \quad p \equiv 1, 2  \quad {\mbox{mod}} \,  4,
  \end{array}
  \right.
  $$
  then $S^+$ and $S^-$ are totally isotropic with respect
  to the symmetric  bilinear form $q(\alpha, \beta)$ defined as the
  top degree component of $\alpha \wedge \sigma (\beta)$ (see \cite{W}).

 A {\it generalized $G_2$-structure} on  a $7$-dimensional manifold $M$ is a reduction from the structure group  $\R^* \times Spin(7,7)$ of
  the bundle $TM  \oplus  T^*M$ to $G_2 \times G_2$.
  Such a structure determines a generalized oriented metric structure
  $(g, b)$, (i.e. a Riemannian  metric $g$, a {\em B}-field $b$ and   an orientation on $V$) and a real scalar function  $\phi$ (the {\em dilaton}).
    Therefore we get a pair of two $G_2$-structures associated with
 two unit spinors $\Psi_{\pm}$ in the irreducible spin representation
  $\Delta = \R^8$ of $Spin(7)$. There is,  up to a scalar,  a unique invariant in $\Lambda^{ev} V^* \otimes \Lambda^{od} V^*$, given by the  box operator
 $$
 \Box_{\rho}: \Lambda^{ev,od} V^*  \to  \Lambda^{od,ev} V^*, \quad \tilde \rho \to e^{\frac{b}{2}} \wedge  \ast_g \sigma (e^{- \frac{b}{2}}  \wedge \tilde \rho).
 $$

 If $\rho$ is a $G_2 \times G_2$-invariant spinor, then  its  {\it companion}   $\hat \rho = \Box_{\rho} \rho$  is still a
 $G_2 \times G_2$-invariant spinor.  To any $G_2 \times G_2$-invariant spinor $\rho$ one can associate a volume  form  $\mathcal Q$ defined by
\begin{equation} \label{volume}
 {\mathcal Q}: \rho \to q (\hat \rho, \rho).
\end{equation}
 Using the isomorphism $\Delta \otimes \Delta \cong \Lambda^{ev,od}$, Witt  in \cite[Proposition 2.4]{W} derived the following  normal form for
 $[\Psi_+ \otimes \Psi_-]^{ev,od}$ in terms of a suitable orthonormal
 basis
 $(e^1, \ldots, e^7)$, namely
 $$
\begin{array}{lcl}
(\Psi_+ \otimes \Psi_-)^{ev} & = & \cos (\theta) + \sin( \theta) (e^{12} + e^{34} + e^{56}) +\\[4pt]
&& \cos  (\theta) (- e^{1367} - e^{1457} - e^{2357} + e^{2467} - e^{1234} - e^{1256} - e^{3456}) +\\[4pt]
&& \sin (\theta )(e^{1357} - e^{1467} - e^{2367} - e^{2457}) - \sin (\theta) e^{123456},\\[4pt]
 (\Psi_+ \otimes \Psi_-)^{odd} &=&  \sin (\theta)  e^7 + \sin (\theta )(-e^{136} - e^{145} - e^{235} + e^{246}) +\\[4pt]
 && \cos  (\theta) (-e^{127} - e^{347} - e^{567} - e^{135} + e^{146} + e^{236} + e^{245}) +\\[4pt]
 && \sin( \theta) (-e^{12347} - e^{12567} - e^{34567}) +\cos (\theta) e^{1234567}, \end{array}$$
 where $\theta$ is the angle between $\Psi_+$ and $\Psi_-$ and $e^{i \ldots j}$ denotes the wedge product $e^i \wedge \ldots \wedge e^j$.

 If  the spinors $\Psi_+$ and $\Psi_-$ are linearly independent, then
 (see Corollary 2.5 of \cite{W})
 $$
\begin{array}{lcl}
(\Psi_+ \otimes \Psi_-)^{ev} & = &  \cos (\theta) + \sin(\theta) \omega  - \cos (\theta) (\psi_-  \wedge \alpha + \frac 12 \omega^2) \\[4pt]
&& + \sin (\theta) \psi_+ \wedge \alpha - \frac{1}{6} \sin (\theta) \omega^3,\\[5pt]
(\Psi_+ \otimes \Psi_-)^{od} & = & \sin (\theta) \alpha -  \cos (\theta) (\psi_+ + \omega \wedge \alpha) - \sin (\theta) \psi_- \\[5pt]
&&- \frac 12 \sin (\theta) \omega^2 \wedge \alpha + \cos  (\theta)
{\mbox {vol}}_g,
\end{array}
 $$
 where $\alpha$  denotes the dual of the unit vector in $V$, stabilized by $SU(3)$,
 $$\omega = e^{12} + e^{34}+ e^{56}$$ is  the fundamental form and
 $\psi_{\pm}$ are  the real and imaginary parts respectively of the complex volume form
 $$
 \psi = (e^1 + i e^2) \wedge (e^3 + i e^4) \wedge (e^5 + i e^6).
 $$
 A $G_2 \times G_2$-invariant spinor $\rho$ is stable in the sense of Hitchin (see \cite{Hi3}), i.e. $\rho$ lies in an open orbit under the
  action of $\R^+  \times Spin(7,7)$.

  By \cite[Theorem 2.9]{W} the generalized $G_2$-structures are in $1-1$ correspondence with lines of spinors
  $\rho$ in $\Lambda^{ev}$ (or $\Lambda^{od} $) whose stabilizer under the action of $Spin(7,7)$ is isomorphic to $G_2 \times G_2$.

 The spinor $\rho$ is called the {\it structure form} of the generalized $G_2$ structure and it can be uniquely written
  (modulo a simultaneous change of sign for $\Psi_+$ and $\Psi_-$) as
 $$\rho = e^{-\phi}  (\Psi_+ \otimes \Psi_-)^{ev}_b ,
 $$
where $b$ is the $B$-field, $\Psi_{\pm} \in \Delta$ are two unit
spinors, the function $\phi$ is  the dilaton and the subscript $b$
denotes the wedge with the exponential $e^{\frac{b}{2} }$.

A {\it {\rm (}topological{\rm )} generalized $G_2$-structure}  over $M$ is a topological $G_2 \times G_2$-reduction of the $SO(7,7)$-principal
 bundle associated with $TM \oplus T^* M$ and it  is characterized by a stable even or odd  spinor $\rho$ which can be viewed as a form.
 This is equivalent to say  that there exists an $SO(7)$-principal fibre bundle which has two $G_2$-subbundles (or equivalently two $G_2{^\pm}$-structures).

 In the sequel we will omit  topological when we will refer to a generalized $G_2$-structure.

Let $H$ be a $3$-form and $\lambda$ be a real, non-zero constant. A generalized $G_2$-structure $(M, \rho)$  is called  {\it strongly integrable} with respect to
$H$ if
$$
d_H \rho = 0, \quad d_H \hat \rho =0,
$$
where $d_H \cdot = d \cdot + H \wedge \cdot$ is the twisted operator
of $d$.  By \cite{W} there are no non-trivial  compact examples  with
a strongly integrable generalized  $G_2$- structure with respect to a
closed $3$-form $H$.

If $$d_H \rho = \lambda  \hat \rho,$$ then the generalized
$G_2$-structure is said  to be {\it weakly  integrable} of {\it
even} or {\it odd} type according to the parity of the form
$\rho$. The constant $\lambda$ (called the {\em Killing number})
is the $0$-torsion form of the two underlying $G_2$-structures.
Indeed, by Corollary 4.6 of \cite{W}, there exist two unique
determined linear connections $\nabla^{\pm}$, preserving the two
$G_2^\pm$-structures, with skew-symmetric torsion $\pm T = \frac
12 db + H$. If the structure is of odd type, then
$$
\begin{array}{l}
d \varphi_+ = \frac{12}{7} \lambda * \varphi_+ + \frac 32 d \phi \wedge \varphi_+ - * T_{27}^+,\\[5pt]
d * \varphi_+ = 2 d \phi \wedge * \varphi_+
\end{array}
$$
and
$$
\begin{array}{l}
d \varphi_- = \frac{12}{7} \lambda * \varphi_- + \frac 32 d \phi \wedge \varphi_- - * T_{27}^-,\\[5pt]
d * \varphi_- = 2 d \phi \wedge * \varphi_-,
\end{array}
$$
where $T_{27}^\pm$ denotes the component of $T$ into the $27$-dimensional irreducible $G_2^\pm$-module
$${\Lambda^3_{27}}^\pm = \{ \gamma \in \Lambda^3 \, \vert \, \gamma \wedge \varphi_+ = \gamma \wedge \varphi_- =0 \}.$$
This is equivalent to say that $e^{-\phi} [\Psi_+ \otimes \Psi_-]$ satisfies the generalized Killing and dilatino equation (see \cite{W,GMPW}).

 In both cases there is a characterization in terms  of the two metric
 connections $\nabla^{\pm}$ with skew symmetric torsion $\pm T$ (see
 \cite[Theorem 4.3]{W}). Indeed, a generalized $G_2$-manifold  $(M,
 \rho)$    is weakly integrable with respect to $H$ if and only if
 $$
  \nabla^{LC} \Psi_{\pm} \pm \frac 14 ( X \lrcorner T) \cdot \Psi_{\pm} =0,\\
  $$
 where $\nabla^{LC}$ is the Levi-Civita connection,  $X \lrcorner$ denotes the contraction by $X$ and the following
 additional conditions are satisfied
 $$
 \left (d \phi \pm \frac 12 ( X \lrcorner T)   \pm \lambda \right) \cdot \Psi_{\pm} =0, \quad
 $$
  if $\rho$ is of even type
  or
 $$
 \left (d \phi \pm \frac 12 ( X \lrcorner T)  + \lambda \right) \cdot \Psi_{\pm} =0, \quad
 $$
  if $\rho$ is of odd type.
  Taking $\lambda = 0$ above equations yield strong integrability  with respect to $H$, instead.

Examples  of  generalized $G_2$-structures  are given by  the {\em
straight} generalized $G_2$-structures, i.e. structures defined by
one spinor $\Psi = \Psi_+ = \Psi_-$. These structures are induced by
a classical $G_2$-structure $(M, \varphi)$ and  are strongly
integrable with respect to a closed  $3$-form $T$ only if the
holonomy  of the metric  associated with $\varphi$ is contained in $G_2$.

If $H$ is closed, then it has to be noted that, in the compact case, the structure form  $\rho$ of  a strongly integrable generalized $G_2$-structure corresponds to a critical point of a functional on stable forms. Indeed, since stability is an open condition, if $M$ is compact  then one can consider the functional
$$
V (\rho) = \int_M {\mathcal Q}(\rho),
$$
where $\mathcal Q$ is defined as in \eqref{volume}.  By \cite[Theorem 4.1]{W} a $d_H$-closed stable form
$\rho$ is a critical point in its cohomology class if and only if $d_H \hat \rho =0$.

Again in the compact case  a $d_H$-exact form $\hat \rho \in \Lambda^{ev,od} (M)$  is  a  critical point  of the functional
 $V$ under some  constraint  if and only if
 $d_H \rho =  \lambda \hat \rho$, for a real non zero constant $\lambda$.

\section{Compact examples of weakly integrable manifolds}

In this section we will construct   examples of compact manifolds endowed with a weakly integrable generalized $G_2$-structure with respect to a closed  $3$-form $H$.

Consider the $7$-dimensional Lie algebra $\mathfrak g$ with structure equations:
$$
\left \{ \begin{array} {l}
d e^1 = a e^{46},\\[3pt]
d e^2 = - \frac 12 a e^{36} - \frac 12 a e^{45} +  \frac 12 a e^{17},\\[3pt]
d e^3 = - \frac 12 a e^{15} +\frac 12 a e^{26} -  \frac 12 a e^{47},\\[3pt]
d e^4 = -a e^{16},\\[3pt]
d e^5 =  \frac 12 a e^{13} -\frac 12 a e^{24} -  \frac 12 a e^{67},\\[3pt]
d e^6 = a e^{14},\\[3pt]
d e^7 =  -\frac 12 a e^{12} -\frac 12 a e^{34} -  \frac 12 a e^{56},
\end{array}
\right.
$$
where $a$ is a real parameter different from zero.

It can be easily checked  that the Lie algebra $\mathfrak g$ is not solvable since $[{\mathfrak g},{\mathfrak g}]  ={\mathfrak g}$ and that  it is unimodular. We can also view  $\mathfrak g$ as the semidirect sum
$$
{\mathfrak g} = {\mathfrak {su}} (2) \oplus_{\delta} \R^4,
$$
where
$$
 {\mathfrak {su}} (2)  = {\mbox {span}} <e_1, e_4, e_6>, \quad   \R^4 = {\mbox {span}} <e_2, e_3, e_5, e_7>
 $$
 and $\delta:   {\mathfrak {su}} (2)  \to {\mathfrak {Der}} (\R^4)$ is given by
 $$
\delta(e_1) = ad_{e_1} = \left( \begin{array}{cccc}   0&0&0&-\frac 12 a\\
 0&0&\frac 12 a&0\\
0&  - \frac 12 a&0&0\\
  \frac 12 a &0&0&0 \end{array} \right),
  $$
    \vskip 0.2cm
 $$
  \delta(e_4) = ad_{e_4} = \left( \begin{array}{cccc}   0&0&\frac 12 a&0\\
  0&0&0&\frac 12 a\\
  - \frac 12 a&0&0&0\\
  0&-  \frac 12 a &0&0\end{array} \right),
  $$
  \vskip 0.2cm
 $$
 \delta(e_6) =  ad_{e_6} = \left( \begin{array}{cccc}   0&-\frac 12 a&0&0\\
  \frac 12 a&0&0&0\\
  0&0&0& \frac 12 a\\
  0&0&- \frac 12 a &0 \end{array} \right).
  $$
 If we identify  $\R^4$ with  the space $\H$  of quaternions, then
 $$
 ad_{e_1} = \frac 12 a L_k, \quad ad_{e_4} = \frac 12 a L_{-j}, \quad ad_{e_6} = \frac 12 a L_{i},
 $$
 where $L_q$ denotes the left multiplication by the quaternion $q$.

Therefore,  the product on the corresponding Lie group $G = SU(2)
\ltimes \H$, for $a = 2$,  is given by
$$
(A, q) \cdot (A', q') = (A A', Aq' + q),  \quad A,A' \in SU(2),
\quad q,q' \in \H,
$$
where  we identify $SU(2)$ with the group of quaternions of unit
norm.

\begin{theorem} \label{example} The Lie group $G = SU(2) \ltimes \H$ admits  compact quotients  $M_{\beta} = G/  \Gamma_{\beta}$, with
 $e^{i \beta}$ primitive p-th root of unity $(p$ prime$)$,  and  $M_{\beta}$ has an invariant weakly integrable generalized $G_2$-structure with respect to a closed $3$-form $H$.
\end{theorem}

\begin{proof}
Consider the discrete subgroup $\Gamma_{\beta} = <A_{\beta} >\ltimes \Z^4$, where $<A_{\beta} >$ is the subgroup of $SU(2)$ generated by
$$
A_{\beta} = \left( \begin{array}{cc} e^{i \beta} &0\\ 0&e^{-i \beta}  \end{array} \right),
$$
with $e^{i \beta}$ primitive p-th root of unity and $p$ prime.

 Then one can check that $\Gamma_\beta$ is a closed subgroup of $G$.  Let $(A', q')$ be any point of
$G$. Thus  $$ [ (A', q')] = \{ (A_{\beta}^m A', A_{\beta}^m q' +
r),  \, m \in \Z \, , r \in \Z^4 \}
$$
is the equivalence class of $(A', q')$. In particular, $[ (A',
q')] = [(A', q' + r)]$ and therefore the restriction of the
projection $\pi: G \to G/\Gamma_{\beta}$  to  $SU(2) \times
[0,1]^4$ is surjective.\newline Then the quotient $M_{\beta} =
(SU(2) \ltimes \H )/ \Gamma_{\beta}$ is a compact manifold.

Consider the invariant metric $g$ on $M_{\beta}$ such that the basis $(e^1, \ldots, e^7)$ is orthonormal
and take the  generalized $G_2$ structure defined by the structure form of odd type
$$
\rho = e^7 - e^{136}- e^{145} - e^{235}+ e^{246} -  e^{12347} - e^{12567} -  e^{34567},
$$
in terms of the basis $(e^1, \ldots, e^7)$.
The
companion  of $\rho$ is
$$
\hat \rho = e^{12} + e^{34} + e^{56}+  e^{1357} -e^{1467} - e^{2367}- e^{2457} - e^{123456}.
$$
Then  the structure form $\rho$ defines a weakly integrable
generalized $G_2$-structure with   respect to  a closed  $3$-form
$H$, i.e. $d_H \rho = \lambda \hat \rho$ ($\lambda$ non-zero
constant),
  if and only if
\begin{equation} \label{weakeq}
\left\{
 \begin{array}{l}
 d e^7 = \lambda \omega,\\[5pt]
 d \psi_- = (H - \lambda \psi_+) \wedge e^7,\\[5pt]
 H \wedge \psi_- = - \frac 13 \lambda \omega^3,
 \end{array}
 \right.
 \end{equation}
 where $\omega, \psi_{\pm}$ are given by
\begin{equation} \label{definitionforms}
\left\{
\begin{array}{lcl}
\omega &= & e^{12} + e^{34} + e^{56},\\[5pt]
\psi_+ &=& e^{135} - e^{146} - e^{236} - e^{245},\\[5pt]
\psi_- &=&  e^{136} + e^{145} + e^{235} - e^{246}.
\end{array}
\right.
\end{equation}
The equations \eqref{weakeq} are satisfied with  $\lambda = - \frac 12 a$ and
$$
 H=     - a e^{146}.
$$
\end{proof}

Observe that $H$ is also co-closed, i.e. $d*H =0$. Moreover, if $a
\leq 1$, $H$ is a calibration in the sense of \cite{HL}.

 In this way we get  compact examples with a weakly
integrable  generalized  $G_2$-structure with respect to the closed
$3$-form $H$.  The induced invariant  metric on $M_{\beta}$ is not flat, since the inner product
$$g = \sum_{i = 1}^7 (e^i)^2$$
on the Lie algebra  $\mathfrak g$  is not flat. Indeed,  the  Ricci
tensor of $g$ is diagonal with respect to the orthonormal basis $(e_1,
\ldots, e_7)$ and its   non zero components are  given by:
$$ Ric (e_1, e_1) = \frac 12 a^2 =Ric (e_4, e_4) = Ric (e_6, e_6).
$$

\section {Link with $SU(3)$-structures in dimension $7$ and
evolution equations}

In this section we will relate the weakly integrable generalized
$G_2$-structures constructed in the previous section with
$SU(3)$-structures in dimension $7$.

 Since the $1$-form $\eta = e^7$ is a contact
form on the Lie algebra ${\mathfrak g}$, then $M_{\beta}$ is   a
contact metric manifold. Moreover, by  \eqref{weakeq} $M_{\beta}$
has an $SU(3)$-structure defined by $(\omega, \eta, \psi = \psi_+
+ i \psi_-)$ such that
\begin{equation} \label{SU3hypo}
\left \{ \begin{array} {l}
d \omega =0,\\[5pt]
d (\psi_\pm \wedge \eta) =0.
\end{array} \right.
\end{equation}
Here we identify $SU(3)$ as the subgroup $SU(3) \times \{1 \}$ of
$SO(7)$.

\smallskip

Note that the $SU(3)$-structures $(\omega, \eta, \psi = \psi_+ + i
\psi_-)$ on $7$-dimensional manifolds for which $d \omega =0$ and
$d (\psi_{\pm}) =0$ where considered in \cite{TV}. In this case
one cannot find any closed $3$-form $H$ such that conditions
\eqref{SU3hypo} are satisfied since $H$ has to be equal to
$\lambda \psi_+$ and the third equation cannot hold. It would be
interesting to investigate if there are other $7$-dimensional
examples endowed with an $SU(3)$-structures which satisfy the
conditions \eqref{SU3hypo} and  giving rise to a weakly integrable
$G_2$-structure with respect to a closed $3$-form $H$.

 In general, let ${\iota}: M^7 \to N^8$ be an embedding
of a an oriented $7$-manifold $M^7$ into a $8$-manifold $N^8$ with
unit normal vector $V$. Then an $SU(4)$-structure  $(\tilde
\omega, \tilde g, \tilde \psi)$ (or equivalently a special almost
Hermitian structure, see e.g. \cite{Ca2}), where  $(\tilde \omega,
\tilde g)$ is a $U(4)$-structure and $\tilde \psi= \tilde \psi_+ +
i \tilde \psi_-$ is  complex  $4$-form of unit norm, defines in a
natural way an $SU(3)$-structure  $(\omega, \eta, g, \psi = \psi_+
+ i \psi_-)$ on $M^7$ given by:
$$
\eta = - V  \lrcorner  \tilde\omega, \quad
\omega = {\iota}^* \tilde \omega, \quad
g = {\iota}^* g, \quad
\psi_+ =  - V  \lrcorner   \tilde\psi_+, \quad
\psi_- =   V \lrcorner \tilde\psi_- .
$$
Then, if $\gamma$ denotes the $1$-form dual to $V$, then  we have
$$
\begin{array} {l}
\tilde \omega =  \omega +  \eta \wedge \gamma,\\[5pt]
\tilde \psi = (\psi_+ + i \psi_-) \wedge (\eta + i \gamma).
\end{array}
$$
The integrability of the $SU(4)$-structure $(\tilde \omega, \tilde
g, \tilde \psi)$ implies  conditions \eqref{SU3hypo}, which can be
viewed as the analogous of the equations defining the hypo
$SU(2)$-structures in dimension 5 (see \cite{CS}).

Vice versa, given an $SU(3)$-structure $(\omega,  \eta,   \psi )$ on $M^7$, an $SU(4)$-structure on $M^7 \times \R$ is defined by
\begin{equation} \label{SU4}
\begin{array} {l}
\tilde\omega = \omega + \eta \wedge dt,\\[5pt]
\tilde \psi = \psi \wedge (\eta + i dt),
\end{array}
\end{equation}
where $t$ is a coordinate on $\R$.

If the $SU(3)$-structure $(\omega, \eta,   \psi )$ on $M^7$  belongs to a one-parameter family of $SU(3)$-structures   $(\omega(t), \eta(t),   \psi(t) )$ satisfying the equations \eqref{SU3hypo} and such that
\begin{equation} \label{evolutions}
\left \{ \begin{array} {l}
\partial_t  \omega(t) = - \hat d \eta(t),\\[5pt]
\partial_t (\psi_+(t) \wedge \eta (t)) = \hat d \psi_-(t),\\[5pt]
\partial_t (\psi_-(t) \wedge \eta (t)) =  -\hat d \psi_+(t),
\end{array} \right.
\end{equation}
for all $t \in (b, c)$, where $\partial_t$ denotes the derivative
with respect to $t$ and $\hat d$ is the exterior differential on
$M^7$, then the $SU(4)$-structure given  by \eqref{SU4} on $M^7
\times (b, c)$  is integrable, i.e. $\tilde \omega$ and $\tilde
\psi$ are both closed. In particular, the associated Riemannian
metric  on $M^7 \times (b, c)$ has holonomy contained in $SU(4)$
and consequently it is Ricci-flat.

For the manifolds $M_{\beta}$ a solution of the evolution
equations \eqref{evolutions} is given by
$$
\begin{array}{l}
\omega(t) = u(t) v(t) (e^{12} + e^{34} + e^{56}),\\[5pt]
\psi_+(t) = u(t) v(t)^2 (e^{135} - e^{236} - e^{245}) - u(t)^3
e^{146},\\[5pt]
\psi_-(t) = u(t)^2 v(t) (e^{136} + e^{145} - e^{246}) + v(t)^3
e^{235},\\[5pt]
\eta(t) = \frac {1}{v(t)^3} e^7,
\end{array}
$$
where $u(t), v(t)$ solve the system of ordinary differential
equations
$$
\left\{ \begin{array}{l} \displaystyle\frac {d}{dt} (u(t) v(t)) =
\frac 12 a
\displaystyle\frac {1}{v(t)^3},\\[10pt]
\displaystyle\frac {d}{dt} \left( \displaystyle\frac {u(t)}{v(t)}
\right) = \frac 12 a v(t)^3\,,
\end{array} \right.
$$
such that $u(0) = v(0) = 1$. The previous system is equivalent to
\begin{equation} \label{odiffsystem}
\left\{ \begin{array}{l}  u'(t) = \frac 14 a
\left(\displaystyle\frac {1} {v(t)^4} +
v(t)^4\right),\\[10pt]
v'(t) = \frac 14 a \left(\displaystyle\frac {1} {u(t) v(t)^3} -
\displaystyle\frac{v(t)^5}{u(t)}\right).
\end{array}
\right.
\end{equation}

Then, by the theorem on existence of solutions for a  system of
ordinary differential equations, one can show that on a open
interval $(b, c)$ containing $t =0$ the system \eqref
{odiffsystem} admits a unique solution $(u(t), v(t))$ satisfying
the initial condition $u(0) = v(0) = 1$. Actually, the solution is given by
$$
u(t) = 1 + \frac 12 a t, \quad v(t) =1.
$$

Hence, we can prove   the following

\begin{theorem}
On the product of $M_{\beta}$ with some open interval $(b, c)$
there exists a Riemannian metric with  discrete holonomy contained in
$SU(4)$.
\end{theorem}

\begin{proof}
The basis of $1$-forms on the manifold $M_{\beta} \times (b,c)$ given by
$$
\begin{array}{l}
E^1 = (1 + \frac 12 a t) e^1, \, \, E^2 =  e^2,  \,\, E^3 = (1 + \frac 12 a t) e^3,\,\, E^4 = (1 + \frac 12 a t) e^4, \\[10pt]
 E^5 =  e^5,\,\, E^6 = (1 + \frac 12 a t) e^6,\,\, E^7 =  e^7,\,\, E^8 = dt
 \end{array}
$$
is orthonormal with respect to the Riemannian metric with holonomy   contained in $SU(4)$. By a direct computation
we have that the non zero Levi-Civita connection 1-forms are given by
$$
\begin{array} {l}
\theta^1_4 =  -\theta^2_3 = \theta^5_7 = \theta^6_8 = \displaystyle\frac {a} {2 + at} E^6, \\[12pt]
  \theta^1_6 = -\theta^2_5  = - \theta^3_7=  -\theta^4_8 =  \displaystyle-\frac {a} {2 + at} E^4,\\[12pt]
    \theta^1_8  =  - \theta^2_7=   \theta^3_5 = \theta^4_6 =  \displaystyle\frac {a} {2 + at} E^1.
  \end{array}
   $$
Therefore, all  the curvature forms $\Omega^i_j$ vanish and consequently the holonomy  algebra is trivial.
\end{proof}

\section{Strong integrability and $SU(3)$-structures in dimension 6}
In this section we are going to consider the structure form $\rho$ of
even type
\begin{equation} \label{rho}
\rho = \omega + \psi_+ \wedge\alpha-\frac{1}{6}\omega^3
\end{equation}
on the product of a $6$-dimensional manifold $N$ endowed with  an $SU(3)$-structure cross $S^1$. We will investigate which type of
$SU(3)$-structures give rise to a strongly integrable generalized
$G_2$-structure with respect to a non-zero $3$-form.

Let  $N$ be a $6$-dimensional manifold. An {\it $SU(3)$-structure} on $N$ is determined by a Riemannian metric $g$, an orthogonal almost complex structure $J$ and a choice of a complex volume form  $\psi = \psi_+ + i \psi_-$  of unit norm.  We will denote by $(\omega, \psi)$  an $SU(3)$-structure, where $\omega$ is the fundamental form defined by
$$
\omega(X, Y) = g (J X, Y),
$$
for any pair of vector fields $X, Y$ on $N$.  Locally one may choose
an   orthornormal basis $(e^1, \ldots, e^6)$ of  the vector cotangent
space $T^*$  such that $\omega$ and $\psi_{\pm}$ are given by
\eqref{definitionforms}.

These forms satisfy the following compatibility relations
$$
\omega \wedge \psi_{\pm} =0, \quad \psi_+ \wedge \psi_- = \frac 23 \omega^3.
$$
The intrinsic torsion of the $SU(3)$-structure belongs to the space (see \cite{CS})
$$
T^* \otimes {\mathfrak {su}} (3)^{\perp} =  {\mathcal W}_1 \oplus  {\mathcal W}_2 \oplus  {\mathcal W}_3 \oplus  {\mathcal W}_4  \oplus  {\mathcal W}_5,
$$
$ {\mathfrak {su}} (3)^{\perp} $ being the orthogonal complement of $ {\mathfrak {su}} (3)$ in  ${\mathfrak {so}} (6)$ and
$$
\begin{array}{ll}
 {\mathcal W}_1 =  {\mathcal W}^+_1 \oplus  {\mathcal W}^-_1, &\quad  {\mathcal W}^{\pm}_1 \cong \R,\\[5pt]
  {\mathcal W}_2 =  {\mathcal W}^+_2 \oplus  {\mathcal W}^-_2, &\quad  {\mathcal W}^{\pm}_2 \cong {\mathfrak {su}}(3),\\[5pt]
  {\mathcal W}_3 \cong [\![{\mathrm S}^{2,0}]\!],  &\quad {\mathcal W}_4  \cong {\mathcal W}_5 \cong T^*,
  \end{array}
  $$
  where $[\![{\mathrm S}^{2,0}]\!]$ denotes the real representation associated with  the space ${\mathrm S}^{2,0}$ of  complex symmetric tensors of type $(2,0)$.\newline
  The components  of the intrinsic torsion of an $SU(3)$-structure can be expressed by (see e.g. \cite{CS,BV})
\begin{equation} \label{intrinsicforms}
\left\{
\begin{array}{lll}
 d \omega &=&  \nu_0 \, \psi_+  +  \alpha_0 \,  \psi_- +   \nu_1  \wedge \omega+ \nu_3,\\[5pt]
d \psi_+ &=&  \frac23 \alpha_0 \, \omega^2 + \pi_1 \wedge \psi_+ -\pi_2 \wedge \omega,\\[5pt]
d \psi_- &= & - \frac 23 \nu_0 \, \omega^2 +  J \pi_1 \wedge \psi_+  - \sigma_2 \wedge   \omega,\\
\end{array}
\right.
\end{equation}
where $\alpha_0 \in {\mathcal W}_1^+$, $\pi_1 \in {\mathcal W}_5$, $\pi_2 \in {\mathcal W}_2^+$, $\nu_0 \in {\mathcal W}_1^-$,
$ \sigma_2 \in  {\mathcal W}_2^-$, $\nu_1 \in {\mathcal W}_4$, $\nu_3 \in {\mathcal W}_3$.

By definition, an $SU(3)$-structure is called {\it integrable} if the intrinsic torsion vanishes. In this case $\omega$ and $\psi$ are both closed. Therefore, the intrinsic torsion measures the failure of the holonomy group of the Levi-Civita connection of $g$ to reduce to $SU(3)$.

If $(\omega, \psi)$ is in the class ${\mathcal W}_2^+$, then  by using \eqref{intrinsicforms} and taking into account  the conditions $d \omega = d \psi_- =0$,  we get  that  the components $\nu_0,  \alpha_0, \sigma_2, \nu_3, \nu_1, \pi_1$ vanish  and hence
$$d \psi_+ =  -\pi_2 \wedge \omega, $$
with $\pi_2$ belonging to the space
\begin{equation} \label{spacepi2}
\begin{array}{lcl}
{\mathcal W}_2^+ &\cong& \{  \gamma \in  \Lambda^2 \quad  \vert  \quad \gamma \wedge \psi_+=0, \quad  * J \gamma = -
\gamma \wedge \omega \}\\[5pt]
&=& \{  \gamma \in  \Lambda^2 \quad \vert \quad J \gamma = \gamma, \quad \gamma \wedge \omega^2 =0 \}.
\end{array}
\end{equation}

By \cite{BV} the scalar curvature ${\mbox {scal}} (g)$ of the metric $g$ is given by:
$$
{\mbox {scal}}  (g)  = - \frac 12 \vert \pi_2 \vert^2 \, .
$$

\medskip

Let $\alpha$ be  a  closed 1-form on $S^1$.  Consider on  the
product $N  \times S^1$,   the  generalized $G_2$-structure
defined by  the structure form of even type $\rho$ given by
\eqref{rho} with companion
$$
\hat \rho = \alpha - \psi_-  - \frac12  \omega^2 \wedge \alpha.
$$
We have the following
\begin{theorem} \label{stronglyint}
Let $(N, \omega, \psi)$ be a $6$-dimensional manifold endowed with an
 $SU(3)$-structure. The structure form $\rho$, given by \eqref{rho},
 defines  a   strongly integrable generalized $G_2$-structure on $N \times S^1$
 with respect to a  $3$-form $H$ $($ non necessarily closed$)$, i.e. $\rho$ satisfies the conditions
 \begin{equation} \label{dequations} d_H \rho  = d_H  \hat \rho =0 \end{equation}
 if and only if
 $N$ is in the class ${\mathcal W}_2^+$ and $H = \pi_2 \wedge \alpha$.
\end{theorem}

\begin{proof} By \eqref{dequations}  we get
$$
\left\{
\begin{array}{l}
d \omega + d (\psi_+ \wedge  \alpha) - \frac16 d(\omega^3)  +
H \wedge \omega + H \wedge \psi^+ \wedge \alpha =0,\\[5pt]
 d \hat  \rho   + H \wedge \hat \rho = - d \psi_-    - \frac12 d(\omega^2 \wedge \alpha) + H \wedge \alpha - H \wedge \psi_-  =0.
\end{array}
\right.
$$
This is equivalent to say:
\begin{equation} \label{strongconditions}
\left\{
\begin{array}{l}
d \omega =0,\\[5pt]
d (\psi_+ \wedge  \alpha) = - H \wedge \omega,\\[5pt]
H \wedge \psi_+  \wedge \alpha=0 \\[5pt]
d \psi_- =  H \wedge \alpha,\\[5pt]
H \wedge \psi_- =0\,.
\end{array}
\right.
\end{equation}
Hence, in particular
$$
d \psi_-  = 0 , \quad  H \wedge \alpha =0.
$$
It follows that $H =  S \wedge f \alpha$, with $S$ a $2$-form on
$N$ and $f$ a function on $S^1$. Since  $d \omega =  0$, we obtain
$$d \psi^+ \wedge \alpha = - S \wedge \omega \wedge f \alpha,
$$ we have that $f$ has to be a constant $k$ and  $$d \psi_+ = - k S
\wedge \omega,$$
 with $k S = \pi_2$. Since $\pi_2$ is a
$(1,1)$-form, then $\pi_2 \wedge \psi_{\pm} =0$. Therefore,
equations \eqref{strongconditions} are satisfied if and only if
$N$ belongs to the class ${\mathcal W}_2^+$.
\end{proof}

Note that $H$ is closed if and only if $d \pi_2 =0$.

\smallskip

Homogeneous examples of $6$-dimensional manifolds with a
$SU(3)$-structure in the class ${\mathcal W}_2^+$ are given in
[8]. There it was proved that the  $6$-dimensional nilmanifolds
$\Gamma \backslash G$ which  carry an invariant $SU(3)$-structures
in the class ${\mathcal W }_2^+$ are the torus, the $\T^2$-bundle
over $\T^4$ and the $\T^3$-bundle over $\T^3$ associated with the
following nilpotent Lie algebras
$$
\begin{array}{l}
(0,0,0,0,0,0),\\[3pt]
(0,0,0,0,12,13),\\[3pt]
(0,0,0,12,13,23),
\end{array}
$$
where  the notation   $(0,0,0,0, 12,13)$ means that  the dual ${\mathfrak g}^*$ of the Lie algebra ${\mathfrak g}$ has a basis $(e^1, \ldots, e^6)$ such that $d  e^i =0, i =1, \ldots, 4$,
 $d e^5 = e^1 \wedge e^2$ and $d e^6 = e^1 \wedge e^3$.

In \cite{W}  the $\T^2$-bundle over $\T^4$ has been considered and it has been  proved that
it admits a $SU(3)$-structure in the class ${\mathcal W}_2^+$.\newline
By \cite{dBT} the $\T^3$-bundle over $\T^3$ admits a family of $SU(3)$-structures in the class ${\mathcal W}_2^+$ given by
$$
\begin{array}{lcl}
\omega &=& e^{16} + \mu e^{25} + (\mu - 1) e^{34},\\[5pt]
\psi_+ &=&  (1 - \mu) e^{124}  + \mu e^{135}   - \mu (\mu - 1) e^{456}  - e^{236},\\[5pt]
\psi_- &=&- \mu (1- \mu) e^{145}  + (\mu - 1) e^{246} + \mu e^{356} + e^{123},
\end{array}
$$
where $\mu$ is  a real number different from $0$ and $1$. Such  a family of $SU(3)$-structures belongs to  the  class ${\mathcal W}_2^+$ with $$
\pi_2 = \mu^2 e^{25} - (\mu- 1)^2 e^{36} - e^{14},
$$
and  $d \pi_2 \neq 0$.
\newline
Manifolds in the class ${\mathcal W}_2^+$ can be also obtained as
hypersurfaces of $7$-dimensional manifolds with a $G_2$-structure.
The $\T^2$-bundle over $\T^4$ can be also be viewed as a
hypersurface of a $7$-dimensional  manifold with a calibrated
$G_2$-structure, i.e. such that the associated stable $3$-form is
closed. Indeed, if $(M, \varphi)$ is a $7$-dimensional manifold
with a calibrated $G_2$-structure, then any hypersurface $\iota: N
\hookrightarrow M$ with unit normal vector $\nu$ such that the Lie
derivative $L_{\nu} \varphi =0$ admits an $SU(3)$-structure
$(\omega, \psi)$ in the class ${\mathcal W}_2^+$ defined by
$$
\begin{array}{l}
\omega = \nu \lrcorner \varphi,\\[3pt]
\psi_+ =  \nu \lrcorner * \varphi,\\[3pt]
\psi_- = \iota^* \varphi.
\end{array}
$$
For general theory on  an oriented  hypersurface of a $7$-dimensional manifold endowed with a $G_2$-structure see \cite{C}.

If we consider the 7-dimensional nilmanifold associated with the
Lie algebra (see \cite{F})
$$
(0,0,0,-13,-23,0,0)
$$
and the hypersurface which is a maximal integral submanifold of the involutive  distribution defined by the $1$-form $e^6$, then one gets the
 $SU(3)$-structure considered above.

 Another example of hypersurface (non  nilmanifold) can be obtained by the  7-dimensional compact manifold $M = X \times S^1$, where $X$ is the compact
 solvmanifold
  considered by Nakamura (see \cite{N}), associated with the solvable Lie algebra
  $$
  (0, 12 - 45, - 13 + 46, 0,15-24,-16+34,0)
  $$
  and endowed with the $G_2$-structure
  $$
  \varphi = e^{147} + e^{357} - e^{267} + e^{136} + e^{125} + e^{234} - e^{456}.
  $$
  The compact hypersurface, maximal integral submanifold of the involutive  distribution defined by the $1$-form $e^7$, has
 an $SU(3)$-structure in the class ${\mathcal W}_2^+$.

\medskip

We will show that,  if the $SU(3)$-structure is not integrable, then
the $2$-form $\pi_2$ cannot be closed. Indeed,

\begin{prop}\label{strongconstruction} Let $N$ be a $6$-dimensional
  manifold endowed with an $SU(3)$-structure $(\omega, \psi)$ in the class ${\mathcal W}_2^+$.
If    $\pi_2$ is  closed,  then the $SU(3)$-structure is integrable. In particular, the associated Riemannian metric $g$ is Ricci flat.
\end{prop}

\begin{proof} As already remarked, $(\omega, \psi)$ is  in the class ${\mathcal W}_2^+$ if and only if \begin{equation} \label{W_2^+}
d \psi_+ =  -\pi_2 \wedge \omega\,,\quad
d \psi_- =  d \omega = 0,
\end{equation}
with $\pi_2$ satisfying the following relations
$$
\begin{array}{ll}
 \pi_2 \wedge \psi_- =0, \quad & * J \pi_2 = -
\pi_2 \wedge \omega\\[3pt]
J \pi_2 = \pi_2, \quad &\pi_2 \wedge \omega^2 =0.
\end{array}
$$
By our assumption  that $\pi_2$ is  closed,  \eqref{W_2^+}  and the above definition of ${\mathcal W}^+_2$ (see \eqref{spacepi2}) we have
$$
0 = d (\pi_2 \wedge \psi_+) =  \pi_2  \wedge d \psi_+ = \pi_2 \wedge d \psi_+ = \vert \pi_2 \vert^2 * 1 \,.
$$
Then  $\pi_2 =0$ and we get   the result.
\end{proof}

In particular, as a consequence we have that  if  $(N, \omega,
\psi)$ is  $6$-dimensional manifold endowed with a (not
integrable)
 $SU(3)$-structure in the class ${\mathcal W}_2^+$, the  $3$-form $H = \pi_2 \wedge \alpha$ on $N \times S^1$ cannot be closed.

\begin{rem} {\rm It has to be noted that, in view of Proposition \ref{strongconstruction}, for $SU(3)$-manifolds in the
class ${\mathcal W}_2^+$, the two conditions
$$
d\pi_2=0\quad \hbox{\rm and}\quad d\psi_+=0
$$
are equivalent.\newline Furthermore, under the conditions of
Proposition \ref{strongconstruction}, the holonomy group of the
metric on the manifold $N$ can be properly contained in $SU(3)$. Indeed, for example, if one
takes the  $6$-manifold  $N=M^4 \times \mathbb{T}^2$, where $(M^4$, $ \omega_1, \omega_2, \omega_3)$
is an hyper-K\"ahler manifold and
$\mathbb{T}^2$ is a $2$-dimensional torus, then an $SU(3)$-structure is defined by
\begin{eqnarray*}
\omega &=& \omega_1 + e^5 \wedge e^6,\\
\psi_+ &=& \omega_2 \wedge e^5 - \omega_3 \wedge e^6,\\
\psi_- &=& \omega_2 \wedge e^6 + \omega_3 \wedge e^5,
\end{eqnarray*}
where $\{e^5, e^6\}$ is an orthonormal coframe on $\mathbb{T}^2$. Since
$$
d \omega_i =0\,,\,\, i = 1,2,3\,,\quad \quad d e^5 = d e^6 =0\,,
$$
we have
$$
d \omega =0\,,\quad d \psi_{\pm} =0\,.
$$
Therefore, the manifold $N$ endowed with the $SU(3)$-structure defined by $(\omega,\psi)$ belongs to the class ${\mathcal W}_2^+$ and
the holonomy of the associated Riemannian metric is strictly contained in $SU(3)$, since the metric is a product.}
\end{rem}
\begin{rem}  {\rm Consider on $N \times \R$ the generalized $G_2$-structure defined by  the structure form $\rho$ given
by \eqref {rho} and let $H$ be a closed non-zero $3$-form.
If we drop  the condition $d_H \hat \rho =0$, then  the $SU(3)$-structure $(\omega, \psi)$ on $N$ has to be in the class
${\mathcal W _2}^+ \oplus {\mathcal W _2}^- \oplus {\mathcal W _5}$ with
$$
d \psi_+  = \pi_1  \wedge \psi_+  - \pi_2 \wedge \omega = - S \wedge \omega,  \quad dS =0.
$$
 Indeed,
$\rho$ is $d_H$-closed  if and only if
$$
\left\{
\begin{array}{l}
d \omega=0\,,\\[5pt]
d  \psi_+   \wedge \alpha =  - H  \wedge \omega\,,\\[5pt]
H \wedge \psi_+ \wedge \alpha =0\,.
\end{array}
\right.
$$
Setting
$$
H = \tilde H + S \wedge \alpha,
$$
with $\tilde H$  and $S$ a $3$-form  and a $2$-form respectively on $N$, then one gets the equivalent conditions:
$$
\left\{
\begin{array}{l}
d \omega =0,\\[5pt]
d \psi_+ = - S \wedge \omega\,,\\[5pt]
\tilde H \wedge \psi_+ = \tilde H \wedge \omega =0\,,\\[5pt]
d S  = d \tilde H =0\,.
\end{array}
\right.
$$
In terms of the components of the intrinsic torsion one has that $\nu_0, \alpha_0,  \nu_1, \nu_3$ vanish and  $$
d \psi_+ = - S  \wedge \omega.
$$
In contrast with the   case of $SU(3)$-manifolds in the class
${\mathcal W}_2^+$  (see Proposition \ref{strongconstruction}),
$6$-dimensional compact examples of this type  may exist, as
showed by the following

\begin{ex} {\rm Consider the $6$-dimensional nilpotent Lie algebra $\mathfrak l$  with structure equations
$$
(0,0,0,0,0,25)
$$
and  the  $SU(3)$-structure given by
$$
\begin{array}{l}
\omega = e^{12} + e^{34} + e^{56},\\[3pt]
\psi = (e^1 + i e^2) \wedge  (e^3 + i e^4) \wedge  (e^5 + i e^6).
\end{array}
$$
Let  $H$ be  the   closed $3$-form
$$
\begin{array} {lcl}
H &= &- e^{457} + a_1 (e^{124} - e^{456}) + a_2 (e^{125} - e^{345})  - a_3 (e^{134}  - e^{156}) + a_4 e^{135} + \\[3pt]
&& a_5 (e^{145} - e^{235}) + a_6 (e^{145} + e^{246}) + a_7 (e^{234} - e^{256}) + a_8 e^{245},
\end{array}
$$
with $a_i \in \R$, $i = 1, \ldots, 8$.
Then $( \omega, \psi)$ induces a structure form  $\rho$ on  a compact quotient of $L \times  \R$, where $L$ is the simply connected nilpotent Lie group
with Lie algebra $\mathfrak l$, by a uniform discrete subgroup. A straightforward computation shows that  $d_H \rho =0$ . }
\end {ex}

}
\end{rem}

\end{document}